\documentclass[11pt]{article}

\usepackage[T2A]{fontenc}
\usepackage[cp1251]{inputenc}
\usepackage{amsfonts}
\usepackage{eufrak}
\usepackage{amssymb}
\usepackage{amsmath}



\begin{document}

\centerline{\textbf{The Heyde characterization theorem on}}

\centerline{\textbf{some locally compact Abelian groups}}

\bigskip

\centerline{\textbf{G.M. Feldman}}

\bigskip

 \makebox[20mm]{ }\parbox{125mm}{ \small {
By the   Heyde  theorem,  the Gaussian distribution
on the real line is characterized by the symmetry of the conditional
distribution of one linear form of of $n$
independent random variables given another.
When $n=2$ we prove analogues
  of this  theorem in the case when independent
  random variables take values
  in a locally compact Abelian group $X$ and coefficients of the
linear forms are topological automorphisms of $X$.}}

\bigskip


\bigskip


\bigskip

\centerline{\textbf{1.   Introduction }}

\bigskip

The well-known  Heyde  theorem states that the Gaussian distribution
on the real line is characterized by the symmetry of the conditional
distribution of one linear form of
$n$ independent random variables given another (\cite{He},
see also \cite[\S\,13.4.1]{KaLiRa}). If $n=2$,
  this theorem can be formulated as follows.

{\bf Theorem A.}  {\it  Let $\xi_1$  and $\xi_2$ be independent random variables with distributions
       $\mu_1$ and $\mu_2$. Let
$\alpha\ne -1$.
If the conditional distribution of  the linear form  $L_2 = \xi_1 +
\alpha\xi_2$  given $L_1 = \xi_1 + \xi_2$ is symmetric, then
  $\mu_1$ and $\mu_2$ are Gaussian distributions.}

Group analogs of the Heyde theorem in the case when random variables
$\xi_j$ take values in a locally compact Abelian group and coefficients
of linear forms are topological automorphisms of the group were studied in
\cite{Fe2}--\cite{Fe3},  \cite{Fe20bb}--\cite{F},   \cite{My2}--\cite{MyF},
see also
\cite[Chapter VI]{Fe5}. Here we continue  this research.

Let $X$ be a second countable locally compact Abelian group. We will
consider only such groups,
without mentioning it specifically.
Denote by ${\rm Aut}(X)$ the group of topological automorphisms of
the group $X$, and by $I$ the identity
automorphism of a group.
Let $\alpha$ be a topological automorphism of the group
$X$ such that   $K={\rm Ker}(I+\alpha)\ne\{0\}$.
Let $\xi_1$ and
$\xi_2$ be independent identically distributed random variables with
values
in the group
      $K$  and distribution  $\mu$.   It is easy
to see that
the conditional distribution of the linear form
    $L_2 = \xi_1 +
\alpha\xi_2$  given $L_1 = \xi_1 + \xi_2$  is symmetric for  any  $\mu$.
Denote by $G$ the subgroup of $X$ generated by all elements of $X$
of order 2. It is obvious, that if $\xi_1$ and
$\xi_2$ are independent  random variables with
values in $G$, then the conditional distribution of the linear form
    $L_2 = \xi_1 + \alpha\xi_2$  given $L_1 = \xi_1 + \xi_2$  is symmetric for  any topological automorphism $\alpha$ of the group
$X$.
Thus, if a group $X$ contains no elements of order 2,
and in this article we will consider only such groups, and
we want to get a theorem, where a class of
distributions
on a group $X$ is
characterized by the symmetry of the conditional distribution of the
linear form
  $L_2 = \xi_1 + \alpha\xi_2$  given $L_1 = \xi_1 + \xi_2$, then we must assume that the condition
 \begin{equation}\label{d1}
{\rm Ker}(I+\alpha)=\{0\}
\end{equation}
  is fulfilled.

On the real line condition (\ref{d1}) is equivalent to the condition
        $\alpha\ne -1$, and it follows from Theorem A that
        for the real line condition (\ref{d1}) is
not only necessary but also sufficient for characterization of the Gaussian distribution by the symmetry
of the conditional distribution of the linear form
$L_2 = \xi_1 + \alpha\xi_2$  given $L_1 = \xi_1 + \xi_2$.
We note that in all the above-mentioned articles devoted to group analogs of the Heyde theorem, where  the linear forms $L_1 = \xi_1 + \xi_2$ and
  $L_2 = \xi_1 + \alpha\xi_2$ are considered, they assumed that a topological automorphism
 $\alpha$ of a group $X$  satisfies  much stronger  conditions than
  ({\ref{d1}), namely
  \begin{equation}\label{E}
I\pm\alpha\in {\rm Aut(X)}.
\end{equation}

In \S 2 we prove that if a locally compact Abelian group $X$   contains no
non-zero  compact elements, then exactly as on the real line, condition
 (\ref{d1}) is necessary and sufficient for characterization of the
 Gaussian distribution on $X$ by the symmetry
of the conditional distribution of the linear form
$L_2$  given $L_1$. In \S 3 we generalize this result to  locally compact Abelian groups $X$ such that
the subgroup of all compact elements of $X$  is  finite and contains no elements of order    $2$. In \S 4  we assume that
the   characteristic functions of independent random variables
$\xi_1$ and
$\xi_2$ do not vanish and prove that if a locally compact Abelian group contains no
elements of order 2, condition
 (\ref{d1}) is necessary and sufficient for characterization of the
 Gaussian distribution on $X$ by the symmetry
of the conditional distribution of the linear form
$L_2$  given $L_1$.

The structure theory for locally compact Abelian groups and the duality theory
will be used in the work (see \cite{Hewitt-Ross}).
Remind some definitions and agree on notation.

 Denote by  $Y$   the character group of the group  $X$, and by
 $(x,y)$ the value of a character   $y \in Y$ at an element  $x
\in X$. If   $G$ is a closed subgroup of    $X$,  denote by
 $A(Y, G) = {\{y \in Y: (x, y) = 1,  \  x \in G \}}$ its annihilator. If  $\alpha\in{\rm Aut}(X)$,
then the adjoint  automorphism $\tilde\alpha\in{\rm Aut}(Y)$ is defined
as follows
   $(x, \tilde\alpha y) = (\alpha x, y)$ for all
 $x \in X$, $y \in Y$. Note that  $\alpha\in {\rm Aut}(X)$
 if and only if $\tilde\alpha\in {\rm Aut}(Y)$.
 Recall that an element
$x\in X$ is said to be compact, if the smallest closed subgroup
of $X$ containing $x$ is compact.

 Denote by $\mathbb{R}$   the group of real numbers, and by
$\mathbb{T}$ the circle group    (the one dimensional torus).

   Let $f(y)$ be a function on the group    $Y$,   and let $h \in
Y$. Denote by   $\Delta_h$   the finite difference operator
$$
\Delta_h f(y) = f(y + h) - f(y), \quad y \in Y.
$$
Let  $n$ be an integer.
 Denote by $f_n:X \mapsto X$ an  endomorphism of the group $X$
 defined by the formula  $f_nx=nx$, $x\in X$.
Put $X^{(n)}=f_n(X)$.

Denote by ${\rm M}^1(X)$   the convolution semigroup of probability distributions on
   $X$.  Let $\mu\in{\rm M}^1(X)$. Denote by  $$\hat\mu(y) =
\int_{X}(x, y)d \mu(x), \quad y\in Y,$$
the characteristic function of the distribution $\mu$, and denote by
$\sigma(\mu)$ the support of $\mu$.

A distribution  $\gamma\in {\rm M}^1(X)$  is called Gaussian
(see \cite[Chapter IV, \S 6]{Pa}),
if its characteristic function is represented in the form
\begin{equation}\label{f1}
\hat\gamma(y)= (x,y)\exp\{-\varphi(y)\}, \quad y\in Y,
\end{equation}
where $x \in X$, and $\varphi(y)$ is a continuous non-negative function
on the group $Y$
 satisfying the equation
 \begin{equation}\label{f2}
\varphi(u + v) + \varphi(u
- v) = 2[\varphi(u) + \varphi(v)], \quad u,  v \in
Y.
\end{equation}
Observe that degenerate distributions are Gaussian.
 Denote by $\Gamma(X)$ the set of Gaussian distributions on
    $X$.  Denote by $E_x$  the degenerate distribution concentrated at an element
$x\in X$.
If $K$ is a compact subgroup of  $X$,  then denote by  $m_K$ the Haar distribution on
 $K$. Denote by  $I(X)$ the set of idempotent distributions on
$X$, i.e. the set of shifts of Haar distributions $m_K$  of compact
subgroups
  $K$ of the group $X.$ For any
  $\mu \in {\rm M}^1(X)$  define the distribution
  $\bar \mu \in {\rm M}^1(X)$
 by the formula $\bar \mu(B) = \mu(-B)$   for any Borel subset   $B$ of $X$.
Note that  $\hat{\bar{\mu}}(y)=\overline{\hat\mu(y)}$.

\bigskip

\centerline{\textbf{2. The Heyde theorem on locally compact Abelian}}

\centerline{\textbf{ groups containing no non-zero compact elements}}

\bigskip

The main result of this section is the following theorem.

{\bf Theorem 1.}
 {\it Let  $X$ be a locally compact Abelian group containing no non-zero compact elements.
  Let  $\alpha$ be a topological automorphism of
    $X$. Let
  $\xi_1$ and  $\xi_2$ be independent random variables with values in the group
       $X$  and distributions
  $\mu_1$ and  $\mu_2$.
The symmetry of the conditional distribution of the linear form
$L_2 = \xi_1 + \alpha\xi_2$ given  $L_1 = \xi_1 +
\xi_2$ implies that   $\mu_1$ and $\mu_2$ are Gaussian distributions if and only if
  $\alpha$ satisfies condition  $(\ref{d1})$.}

To prove Theorem 1 we need some lemmas.

{\bf  Lemma 1} (\cite[Lemma 16.1]{Fe5}). {\it  Let $X$ be a locally compact Abelian group, $Y$ be
its character group. Let $\alpha$ be a topological automorphism of $X$.
 Let $\xi_1$ and $\xi_2$ be independent random variables with values in the group
       $X$  and distributions
 $\mu_1$ and $\mu_2$.
 The conditional distribution of the linear form
 $L_2 = \xi_1 + \alpha\xi_2$
 given $L_1 = \xi_1 + \xi_2$ is symmetric if and only
 if the characteristic functions
 $\hat\mu_j(y)$ satisfy the equation}
\begin{equation}\label{42}
\hat\mu_1(u+v )\hat\mu_2(u+\tilde\alpha v )=
\hat\mu_1(u-v )\hat\mu_2(u-\tilde\alpha v), \quad u, v \in Y.
\end{equation}

{\bf  Corollary 1.} {\it Let $X$ be a locally compact Abelian group.
 Let $\xi_1$ and $\xi_2$ be independent
 random variables with values in
       $X$ and distributions  $\mu_1$ and $\mu_2$. Then the following statements hold.

$1$. If $\mu_1=\mu_2$, then the conditional distribution of the linear form
 $L_2 = \xi_1 -\xi_2$
 given $L_1 = \xi_1 + \xi_2$ is symmetric.

$2$. If the group $X$ contains no elements of order  $2$ and  the conditional distribution of the linear form
 $L_2 = \xi_1 -\xi_2$
 given $L_1 = \xi_1 + \xi_2$ is symmetric, then $\mu_1=\mu_2$}.

{\bf Proof}. 1. Obviously, if  $\mu_1=\mu_2$  and $\alpha=-I$, then the characteristic functions
$\hat\mu_j(y)$ satisfy equation
 $(\ref{42})$. By Lemma  1 this implies that the conditional distribution of the linear form
 $L_2 = \xi_1 -\xi_2$
 given $L_1 = \xi_1 + \xi_2$ is symmetric.

2. By Lemma  1 the characteristic functions
$\hat\mu_j(y)$ satisfy equation
 \begin{equation}\label{13.02.2016.1}
\hat\mu_1(u+v )\hat\mu_2(u- v )=
\hat\mu_1(u-v )\hat\mu_2(u+ v), \quad u, v \in Y.
\end{equation}
Substituting $u=v=y$ in  (\ref{13.02.2016.1}), we get
\begin{equation}\label{13.02.2016.2}
\hat\mu_1(2y)=\hat\mu_2(2y), \quad y \in Y.
\end{equation}
 Since the group $X$ contains no elements of order 2, the subgroup $Y^{(2)}$ is dense in $Y$. Hence, (\ref{13.02.2016.2}) implies that $\hat\mu_1(y)=\hat\mu_2(y)$, $y\in Y$, and so $\mu_1=\mu_2$.
   $\Box$

 It is convenient for us to formulate as a lemma the following well-known statement
 (see e.g. \cite[Proposition 2.13]{Fe5}).

{\bf Lemma 2}. {\it Let $X$ be a locally compact Abelian group,
$Y$ be its character group. Let    $\mu\in{\rm
M}^1(X)$. Then the set  $E=\{y\in Y:\ \hat\mu(y)=1\}$ is a closed
subgroup of
   $Y$,   and the characteristic function
$\hat\mu(y)$ is $E$-invariant, i.e.  $\hat\mu(y)$ takes a constant value on each coset
of the subgroup $E$ in
  $Y$. Moreover, $\sigma(\mu)\subset A(X,E)$.}

The following statement  is contained in
the proof of Proposition 1 in  \cite{Fe3}, see also \cite[Theorem 17.19]{Fe5}.

{\bf  Lemma 3.} {\it Let $Y$ be a connected compact Abelian group,
$W$ be a neighbourhood of
zero of $Y$. Let $P(y)$ be a continuous function satisfying the equation
\begin{equation}\label{31a}
\Delta_h^kP(y)=0, \quad  y, h\in W,
\end{equation}
 for some natural  $k$,  and the condition  $P(0)=0$.
    Let $\psi(y)$ be a characteristic function on the group
 $Y$ such that the  representation
\begin{equation}\label{31}
\psi(y)=e^{P(y)}, \quad y\in W,
\end{equation}
holds.
   Then $\psi(y)=1,$ $y\in Y$.}

{\bf Lemma 4.} {\it  Let  $X$ be a torsion free discrete Abelian group.
Let  $\alpha$ be an automorphism of
    $X$,  satisfying condition $(\ref{d1})$.
 Let $\xi_1$ and  $\xi_2$ be independent random variables with
 values in the group
 $X$ and distributions       $\mu_1$ and $\mu_2$.
If the conditional distribution of the linear form    $L_2 = \xi_1 +
\alpha\xi_2$  given $L_1 = \xi_1 +\xi_2$ is symmetric,
then $\mu_1$ and $\mu_2$ are degenerate distributions.}

{\bf Proof.} Note that if an automorphism   $\alpha$ of the group $X$ satisfies
 conditions $(\ref{E})$,
then  Lemma 4 is a particular case of a more general statement (\cite[ Proposition  1]{Fe3}, and
see also
  \cite[Theorem 17.19]{Fe5}). We follow the scheme of the proof
 of  Theorem 17.19 in \cite{Fe5}.

 Denote by   $Y$  the character group of the group  $X$. By Lemma  1,
 the symmetry of the conditional distribution of
 the linear form
 $L_2$ given  $L_1$ implies that the characteristic functions
$\hat\mu_j(y)$  satisfy equation $(\ref{42})$. Put
$\nu_j = \mu_j* \bar \mu_j$.
Then  $\hat \nu_j(y) = |\hat \mu_j(y)|^2 \ge 0,$  $y \in Y$, $j=1, 2$,
and the characteristic functions
   $\hat \nu_j(y)$ also satisfy equation
   $(\ref{42})$.  Thus  in the course of the proof of Lemma 4, we may
   suppose without loss
   of generality that   $\hat\mu_j(y)\ge 0$, $j=1, 2$. We will
   prove that in this case  $\mu_1=\mu_2=E_0$.

Let  $U$  be a neighbourhood of zero of the group  $Y$ such
that    $\hat\mu_j(y)> 0$  for
$ y \in U$, $j=1, 2$. Put
  $\varphi_j(y) = - \log \hat\mu_j(y)$, $y \in U$.  Let  $V$ be a symmetric neighbourhood of zero
  of  $Y$ such that for any automorphism    $\lambda_j \in \{I, \tilde\alpha \}$
  the inclusion
 \begin{equation}\label{11}
\sum_{j=1}^8 \lambda_j (V) \subset U
\end{equation}
 holds.
It follows from $(\ref{42})$ that the functions   $\varphi_j(y)$
satisfy the equation
\begin{equation}\label{12}
\varphi_1(u+v) + \varphi_2(u+\tilde\alpha v)- \varphi_1(u-v)
- \varphi_2(u-\tilde\alpha v)=0, \quad u, \ v \in V.
  \end{equation}

We use the finite difference method to solve equation (\ref{12}) (see e.g. \cite[Ch. VI]{Fe5}).
Since (\ref{11}) holds, arguments of the functions
   $\varphi_j(y)$ are in  $U$.

Let $k_1$  be an arbitrary element of  $V$. Put $h_1=\tilde\alpha k_1$, and substitute in
  (\ref{12}) $u+h_1$ for $u$ and $v+k_1$ for $v$. Subtracting equation
  (\ref{12}) from the resulting equation, we get
\begin{equation}\label{13}
\Delta_{l_{11}}\varphi_1(u+v) + \Delta_{l_{12}}\varphi_2(u+\tilde\alpha v)
-\Delta_{l_{13}}\varphi_1(u-v)=0, \quad u, v \in
V,
\end{equation}
where $l_{11}= (I+\tilde\alpha)k_1$, $l_{12}=2 \tilde\alpha k_1$, $l_{13}= (\tilde\alpha-I)k_1.$

Let  $k_2$  be an arbitrary element of  $V$. Put $h_2=k_2$, and substitute
in   (\ref{13}) $u+h_2$  for $u$ and $v+k_2$ for $v$.
Subtracting equation
  (\ref{13}) from the resulting equation, we obtain
\begin{equation}\label{14}
\Delta_{l_{21}}\Delta_{l_{11}}\varphi_1(u+v) +
\Delta_{l_{22}}\Delta_{l_{12}}\varphi_2(u+\tilde\alpha v) =
0,  \quad u, v \in V,
\end{equation}
where $l_{21}=2 k_2$, $l_{22}=(I+\tilde\alpha)k_2.$

Let $k_3$ be an arbitrary element of $V$. Put $h_3 =-\tilde\alpha
k_3$. Substitute in
(\ref{14}) $u+h_3$ for $u$ and $v+k_3$ for $v$.   Subtracting equation
  (\ref{14}) from the resulting equation, we get
\begin{equation}\label{15}
\Delta_{l_{31}}\Delta_{l_{21}}\Delta_{l_{11}}\varphi_1(u+v)
= 0 \quad u, v \in V,
\end{equation}
where $l_{31}= (I-\tilde\alpha)k_3.$
Putting in (\ref{15}) $v=0$,  we find
\begin{equation}\label{16}
\Delta_{l_{31}}\Delta_{l_{21}}\Delta_{l_{11}}\varphi_1(u) = 0,
 \quad u \in V.
\end{equation}

Reasoning similarly we obtain from (\ref{14}) that the function $\varphi_2(y)$
satisfies the equation
\begin{equation}\label{new16}
\Delta_{l_{32}}\Delta_{l_{22}}\Delta_{l_{12}}\varphi_2(u) = 0,
 \quad u \in V,
\end{equation}
where $l_{32}=-(I-\tilde\alpha)k_3$.

Note that    $X$ is a torsion free discrete Abelian group if and only if  $Y$ is a connected compact
 Abelian group (\cite[(23.17), (24.25)]{Hewitt-Ross}). Since
  $\alpha$ satisfies condition
    $(\ref{d1})$, the subgroup
    $(I+\tilde\alpha)(Y)$ is dense in
    $Y$ (\cite[(24.41)]{Hewitt-Ross}), and
    taking into account that   $Y$ is a compact group, we have
    \begin{equation}\label{16n}
(I+\tilde\alpha)(Y)=Y.
\end{equation}
This implies that  $I+\tilde\alpha$ is an open homomorphism
  (\cite[(5.29)]{Hewitt-Ross}), and hence,  $l_{11}(V)$
is a neighbourhood of zero of  $Y$.
Since $Y$ is a connected group,  $Y$ is a divisible group (\cite[(24.25)]{Hewitt-Ross}),
i.e. for any natural   $n$  the equality  $Y^{(n)} = Y$ holds. In
particular,
$Y^{(2)} = Y$.  This implies that $f_2 : Y \mapsto Y$ is an open homomorphism (\cite[(24.41)]{Hewitt-Ross}),
 and hence, $l_{21}(V)$ is also a neighbourhood of zero of the group $Y$.

 Put $L=(I-\tilde\alpha)(Y).$ Then $L$ a connected compact
 subgroup of the group
   $Y$,  and a continuous homomorphism  $I-\tilde\alpha$,  as a homomorphism
   from $Y$ to $L$, is open (\cite[(5.29)]{Hewitt-Ross}). It follows
   from this that  $l_{31}(V)$ is
a neighbourhood of zero of the group $L$. Taking into account $(\ref{16})$,
it follows from what had been said that
 there exists a neighbourhood of zero  $W$  of the group $L$ such that the function $\varphi_1(y)$
 satisfies the equation
 \begin{equation}\label{17}
\Delta_h^3\varphi_1(y)=0, \quad y, h \in W.
\end{equation}
 Considering the group $L$ instead of $Y$ and applying  Lemma 3 to the group
 $L$, we obtain
  from  $(\ref{17})$ that
 $\hat\mu_1(y) = 1$ for  $y\in L$.
Hence, by Lemma 2, $\sigma(\mu_1)\subset A(X, L)$.

Put $G=A(X, L)$.
Reasoning similarly we get from $(\ref{new16})$ that $\sigma(\mu_2)\subset G$.
It follows from $L=(I-\tilde\alpha)(Y)$ that   $G={\rm Ker}(I-\alpha)$ and,
obviously,
$\alpha(G)=G$.
Since  $\sigma(\mu_j)\subset  G$, $j=1, 2$,  the independent random
variables   $\xi_1$ and  $\xi_2$
 take values in  $G$. Denote by $N$ the character group of the group $G$.  The characteristic functions of distributions  $\mu_j$, considering
 as distributions on $G$,  satisfy equation (\ref{42}) on  $N.$ Moreover,  since the restriction of $\alpha$
  to $G$ is the identity automorphism, equation
   (\ref{42}) on  $N$  take the form
\begin{equation}\label{44}
\hat\mu_1(u+v )\hat\mu_2(u+v )=
\hat\mu_1(u-v )\hat\mu_2(u-v), \quad u, v \in N.
\end{equation}
 Substituted in (\ref{44}) $u=v=y$, and taking into account that $\hat\mu_j(y)\ge 0$,
 we obtain
\begin{equation}\label{881}
\hat\mu_1(2y)=\hat\mu_2(2y)=1, \quad y\in N.
\end{equation}

Since $G$ is a torsion free discrete Abelian group,
$N$ is a connected compact Abelian group, and hence,
$N^{(n)}=N$  for any natural  $n$ (\cite[(24.25)]{Hewitt-Ross}). In particular,
  $N^{(2)}=N$.  It follows from  (\ref{881})  that then
$\hat\mu_1(y)=\hat\mu_2(y)=1$  for any  $y\in N$.
 This implies that $\mu_1=\mu_2=E_0$.
$\Box$

 Taking into account that   (\ref{16n}) implies (\ref{d1}) and considering
 Lemma 1,  the following statement follows from Lemma 4.

{\bf Corollary 2.} {\it  Let
 $Y$ be a connected compact Abelian group,   $X$ be
its character group. Let $\tilde\alpha$ be a topological automorphism of
the group   $Y$,  satisfying condition  $(\ref{16n})$.
 Let ${\mu_j\in {\rm M}^1(X)}$, $j=1, 2$, and assume that
 the characteristic functions $\hat\mu_j(u)$ satisfy equation  $(\ref{42})$.
  Then $\hat\mu_j(y)=(x_j, y)$, where  $x_j\in X$, $j=1, 2$.}

 {\bf  Lemma 5} (\cite[Lemma 6]{My1}). {\it Let  $X$ be a locally compact Abelian group. Let $\alpha$
be a topological automorphism of $X$.
Let  $\xi_1$ and  $\xi_2$ be independent random variables with values in
the group $X$.
   If the conditional distribution of the linear form    $L_2 = \xi_1 +
\alpha\xi_2$  given $L_1 = \xi_1 +
\xi_2$  is symmetric,   then the linear forms
$M_1=(I+\alpha)\xi_1+2\alpha\xi_2$ and
$M_2=2\xi_1+(I+\alpha)\xi_2$  are independent.}

We formulate as a lemma the structure theorem for locally compact
Abelian groups.

{\bf Lemma 6} (\cite[(24.30)]{Hewitt-Ross}).  {\it   Any locally compact Abelian
group $X$ is topologically isomorphic to a group of the form  $ \mathbb{R}^n\times G,$
where $n\ge 0$, and the group  $G$  contains a compact open subgroup}.

{\bf  Proof of Theorem 1.} Necessity. Assume that  $K={\rm Ker}(I+\alpha)\ne\{0\}$.
Let $\xi_1$ and $\xi_2$ be independent identically distributed random variables
with values in $K$   and distribution $\mu$. Applying  Corollary  1
to the group $K$, we get that the conditional distribution of
the linear form $L_2 = \xi_1 - \xi_2$ given
  $L_1 = \xi_1 + \xi_2$ is symmetric.
It is obvious that   $\alpha x=-x$  for all
 $x\in K$. Thus,    if we consider the random variables  $\xi_1$
 and  $\xi_2$,  as
the random variables with values in $X$, then the conditional
distribution of the linear form
     $L_2 = \xi_1 +
\alpha\xi_2$  given $L_1 = \xi_1 + \xi_2$  is symmetric.  Taking into
account that
$\mu$ is an arbitrary distribution, the necessity is proved.

Note that we did not use  in this reasoning the fact that the group $X$ contains no
non-zero compact elements.

Sufficiency. Denote by   $Y$  the character group of the group  $X$.
Since   $X$
contains no non-zero compact elements, it follows from Lemma 6 that the group
$X$  is topologically isomorphic to a group of the form
$\mathbb{R}^n\times D$, where
$n\ge 0$,  and $D$ is a torsion free discrete Abelian group.
Denote by
$(t, d)$, where $t\in \mathbb{R}^n$, $d\in D,$ elements of the group $X$.
Denote by $C$ the character group of the group $D$. Then $C$ is a connected compact Abelian group
     (\cite[(23.17), (24.25)]{Hewitt-Ross}),  and the group $Y$
   is topologically isomorphic to a group of the form  $\mathbb{R}^n\times C$.
In order not to complicate the notation  we assume that
$X=\mathbb{R}^n\times D$ and
$Y=\mathbb{R}^n\times C$.

By Lemma  1,   the symmetry of the conditional distribution of the
linear form     $L_2$
 given  $L_1$ implies that the characteristic functions
$\hat\mu_j(y)$ satisfy equation  $(\ref{42})$.
Set $\nu_j = \mu_j* \bar \mu_j$.  Then $\hat \nu_j(y) = |\hat \mu_j(y)|^2 \ge 0,$  $y \in Y$,
and  the characteristic functions   $\hat \nu_j(y)$
 also satisfy equation  $(\ref{42})$.

 Consider the factor-group $X/\mathbb{R}^n$.
Denote by    $[x]$ its elements.
Since   $\mathbb{R}^n$ is the connected component of zero of the group    $X$,
we have $\alpha(\mathbb{R}^n)=\mathbb{R}^n$. Hence,  $\alpha$  induces
  an automorphism $\hat\alpha$ on the factor-group $X/\mathbb{R}^n$ by the formula
 $\hat\alpha [x]=[\alpha x]$. Verify that
 \begin{equation}\label{091}
{\rm Ker}(I+\hat\alpha)=\{0\}.
\end{equation}
 Take $x_0\in X$ such that $[x_0]\in {\rm Ker}(I+\hat\alpha)$. We have
$(I+\hat\alpha)[x_0]=0$, i.e. $[(I+\alpha)x_0]=0$.
This implies that
\begin{equation}\label{Angers3}
(I+\alpha)x_0=x'
\end{equation}
 for some $x'\in \mathbb{R}^n$. It is obvious that
$(I+\alpha)(\mathbb{R}^n)\subset\mathbb{R}^n$. It follows from $(\ref{d1})$ that
the restriction of the continuous endomorphism
$I+\alpha$ of the group $X$  to the
subgroup $\mathbb{R}^n$ is a topological automorphism of the group   $\mathbb{R}^n$.
 Hence, $x'=(I+\alpha) \tilde x$
for some  $\tilde x\in \mathbb{R}^n$.
 In view of (\ref{Angers3}), this implies
$(I+\alpha)x_0=(I+\alpha)\tilde x $. Taking into account
$(\ref{d1})$, it follows from this that  $x_0=\tilde x$.
 Hence,
  $[ x_0]=0$.  So, we proved $(\ref{091})$.

  Since $Y=\mathbb{R}^n\times C$, and  $C$ is a compact group, obviously, $\tilde\alpha(C)=C$.
 Take into account that the character group of the factor-group $X/\mathbb{R}^n$
 is topologically isomorphic in a natural way to the group $A(Y, \mathbb{R}^n)=C$ (\cite[(23.25)]{Hewitt-Ross}), and note that the automorphism adjoint to  $\hat\alpha$ acts on $C$
  as
  the restriction of  $\tilde\alpha$ to $C$.
In view of  $\tilde\alpha(C)=C$,  we can
consider the restriction of equation  $(\ref{42})$  for the characteristic functions
  $\hat \nu_j(y)$ to  $C$. It follows from  $(\ref{091})$ that the subgroup
  $(I+\tilde\alpha)(C)$ is dense in $C$ (\cite[(24.41)]{Hewitt-Ross}).
  Since $C$ is a compact group, this implies that $(I+\tilde\alpha)(C)=C$,
  i.e. the restriction of the topological automorphism $\tilde\alpha$ to $C$
  satisfies
  condition $(\ref{16n})$.
  Apply Corollary 2, considering the group $C$ instead of $Y$.
   Taking into account that $\hat\nu_j(y)\ge 0,$ $j=1, 2,$ we get
   that
$\hat\nu_1(y)=\hat\nu_2(y)=1$ for $y\in C.$ Hence, by Lemma 2,
$\sigma(\nu_j)\subset A(X, C)=\mathbb{R}^n$, $j=1, 2$.

Denote by   $\zeta_1$  and $\zeta_2$   independent random variables with values
in the group $X$  and distributions   $\nu_1$ and $\nu_2$.
  Since the characteristic functions  $\hat \nu_j(y)$ satisfy equation  $(\ref{42})$,
 by Lemma 1,  the conditional distribution of the linear form
 $\tilde L_2 = \zeta_1 +
\alpha\zeta_2$  given $\tilde L_1 = \zeta_1 +
\zeta_2$ is symmetric.
 In view of  $\sigma(\nu_j)\subset \mathbb{R}^n$,
 the independent random variables $\zeta_1$   and $\zeta_2$ takes values in    $\mathbb{R}^n$.
   Taking into account that  $\alpha(\mathbb{R}^n)=\mathbb{R}^n$,
   denote by $\alpha_{_{\mathbb{R}^n}}$  the restriction of $\alpha$ to  $\mathbb{R}^n$.
   Obviously, $\alpha_{_{\mathbb{R}^n}}\in {\rm Aut}(\mathbb{R}^n)$, and the conditional distribution of the
   linear form $L'_2 = \zeta_1 +
\alpha_{_{\mathbb{R}^n}}\zeta_2$ given  $L'_1 = \zeta_1 +
\zeta_2$ is symmetric.  In view of Lemma 5, this implies that
 the linear forms
$M_1=(I+\alpha_{_{\mathbb{R}^n}})\zeta_1+2\alpha_{_{\mathbb{R}^n}}\zeta_2$  and
$M_2=2\zeta_1+(I+\alpha_{_{\mathbb{R}^n}})\zeta_2$ are independent.

 As has been noted above in view of $(\ref{d1})$,   the restriction of the continuous endomorphism
  $I+\alpha$ of the group $X$ to the subgroup  $\mathbb{R}^n$ is a
  topological automorphism of
  $\mathbb{R}^n$, i.e. ${I+\alpha_{_{\mathbb{R}^n}}\in {\rm Aut}(\mathbb{R}^n)}$. Hence,
all coefficients of the linear forms $M_1$ and $M_2$   are topological automorphisms of
the group
 $\mathbb{R}^n$. Applying the well-known Ghurye--Olkin theorem about characterisation
 of the Gaussian vectors  in the space $\mathbb{R}^n$ by
 the independence of linear
 forms of independent random vectors (\cite{GhurO}, see also \cite[Theorem 3.2.1]{KaLiRa}),
 we get that $\nu_j\in\Gamma(\mathbb{R}^n)$, $j=1, 2.$

It is easy to see that we can change the distributions $\mu_j$ by their shifts $\mu'_j$
in such a way that $\mu'_j*\bar\mu'_j=\nu_j$, and $\sigma(\mu'_j)\in \mathbb{R}^n$, $j=1, 2$.
 It follows from the Cramer theorem about decomposition of the Gaussian distribution
 in the space $\mathbb{R}^n$  that $\mu'_j\in \Gamma(\mathbb{R}^n)$ (\cite[Theorem 6.3.2]{LiO}).
 Hence, $\mu_j\in \Gamma(X)$, $j=1, 2$.
$\Box$

{\bf  Corollary 3.} {\it Let  $X$ be a locally compact Abelian group containing no non-zero compact elements.
Let $\alpha_j, \beta_j\in {\rm Aut}(X)$.
Let $\xi_1$ and  $\xi_2$ be independent random variables
with values in the group $X$  and distributions
  $\mu_1$ and  $\mu_2$.
The symmetry of the conditional distribution of the linear form
$L_2=\beta_1\xi_1+\beta_2\xi_2$ given $L_1=\alpha_1\xi_1+\alpha_2\xi_2$  implies that $\mu_1$ and $\mu_2$ are Gaussian distributions
if and only if}
\begin{equation}\label{Kh1}
{\rm Ker} (I+\alpha_1\beta_1^{-1}\beta_2\alpha_2^{-1})=\{0\}.
\end{equation}

{\bf Proof}. Necessity. Set $K={\rm Ker} (I+\alpha_1\beta_1^{-1}\beta_2\alpha_2^{-1})$,
 and assume that $K\ne \{0\}$. Let $\eta_1$ and
$\eta_2$ be independent identically distributed random variables with values in the group
$K$  and distribution  $\mu$. Put $\xi_j=\alpha_j^{-1}\eta_j$, $j=1, 2$, and
check that the conditional distribution of the linear form
$L_2=\beta_1\xi_1+\beta_2\xi_2$ given $L_1=\alpha_1\xi_1+\alpha_2\xi_2$ is symmetric.
This statement is equivalent to the statement that
the conditional distribution of the linear form
$\tilde L_2=\beta_1\alpha_1^{-1}\eta_1+\beta_2\alpha_2^{-1}\eta_2$ given $\tilde
L_1=\eta_1+\eta_2$ is symmetric. Taking into account  that
$\beta_2\alpha_2^{-1}x=-\beta_1\alpha_1^{-1}x$  for any $x\in K$, the linear form
$\tilde L_2$ takes the form $\tilde L_2=\beta_1\alpha_1^{-1}(\eta_1-\eta_2)$,
and the symmetry of the conditional distribution of $\tilde L_2$ given $\tilde
L_1$ follows from Corollary 1. Since we can consider $\eta_1$ and
$\eta_2$ as random variables with values in $X$, and $\mu$ is an arbitrary distribution,
the necessity of condition (\ref{Kh1}) is proved.

To prove sufficiently set $\eta_j=\alpha_j\xi_j$, $j=1, 2$.
The symmetry of the conditional distribution of the linear form
$L_2=\beta_1\xi_1+\beta_2\xi_2$ given $L_1=\alpha_1\xi_1+\alpha_2\xi_2$  implies that
the conditional distribution of the linear form
$\tilde L_2=\eta_1+\alpha_1\beta_1^{-1}\beta_2\alpha_2^{-1}\eta_2$ given
$\tilde L_1=\eta_1+\eta_2$ is symmetric. By Theorem 1, the random variables
$\eta_1$ and $\eta_2$ have have Gaussian distributions. Hence, $\xi_1$ and $\xi_2$
also have Gaussian distributions. Sufficiency of condition (\ref{Kh1}) is also
proved. $\Box$

{\bf  Remark 1.}  It is easy to see that in the class of locally compact Abelian groups
containing no non-zero compact elements, in particular, in the class of torsion
free discrete Abelian groups,
 generally speaking, $(\ref{d1})$ does not imply that $I+\alpha\in {\rm Aut}(X)$. Indeed, let
$$
X=\left\{{m\over {2^n}}: \ n=0, 1, \dots; \ m=0, \pm 1, \pm 2, \dots\right\}
$$
 be a subgroup of the additive group of the rational numbers considering in the discrete topology.
  Let $\alpha$ be an automorphism of the group
 $X$ of the form
 $\alpha x = 2x$, $x\in X$. It is obvious that  $I+\alpha\not\in {\rm Aut}(X)$, whereas
 condition $(\ref{d1})$ holds.

\bigskip

\centerline{\textbf{3.  The Heyde theorem on locally compact Abelian groups}}

\centerline{\textbf{with a finite
subgroup of all compact elements}}

\bigskip

 Let $X$ be a locally compact Abelian group.
 It is well known that the support of any Gaussian distribution on the group
 $X$ is a coset of a connected subgroup  in $X$ (\cite[Chapter IV]{Pa}).
So,  if $X$  is a totally disconnected group, then Gaussian distributions on  $X$
 are degenerate. For such groups the role of  Gaussian distributions play idempotent
 distributions. Taking this into account, the following statement can be considered as the Heyde theorem for finite Abelian groups
(compare with Theorem A).

{\bf Theorem  B} (\cite{Fe2}, see also   \cite[Theorem 17.1 and Remark 17.5]{Fe5}).
  {\it  Let $X$ be a finite Abelian group containing no elements of order
   $2$. Let $\alpha$ be an automorphism of the group
$X$ such that  $I+\alpha\in {\rm Aut(X)}$. Let  $\xi_1$ and  $\xi_2$ be independent random variables with values in
       $X$   and distributions
  $\mu_1$ and $\mu_2$.
 If the conditional distribution of the linear form  $L_2 = \xi_1 + \alpha\xi_2$
given  $L_1 = \xi_1 + \xi_2$ is symmetric,  then $\mu_j\in I(X)$, $j=1, 2$.}

Note that for a finite Abelian group $X$ the condition  $I+\alpha\in {\rm Aut(X)}$
is equivalent to   condition $(\ref{d1})$.
We prove in this section the following
generalization of Theorem 1 and Theorem B.

{\bf Theorem 2.} {\it Let  $X$ be a locally compact Abelian group, containing no elements of order    $2$. Let $F$ be the subgroup of all compact elements of $X$. Assume that
$F$   is a finite group.
Let  $\alpha$ be a topological automorphism of
    $X$. Let
  $\xi_1$ and  $\xi_2$ be independent random variables with values in the group
       $X$  and distributions
  $\mu_1$ and  $\mu_2$.
The symmetry of the conditional distribution of the linear form
$L_2 = \xi_1 + \alpha\xi_2$ given  $L_1 = \xi_1 +
\xi_2$ implies that   $\mu_j\in\Gamma(X)*I(X)$, $j=1, 2$, if and only if
  $\alpha$ satisfies condition  $(\ref{d1})$.}

To prove Theorem 2 we need the following lemma.

{\bf Lemma 7} (\cite{Fe14}), see also \cite[Theorem 13.3]{Fe5}). \textit{
  Let $X=\mathbb{R}^n\times F$, where $F$
is a finite Abelian group. Let $\alpha_j, \beta_j$ be
topological automorphisms of  $X$.  Let   $\xi_1$ and $\xi_2$ be
independent random variables with values in the group  $X$ and
distributions $\mu_1$ and $\mu_2$.
If the linear forms
$L_1 = \alpha_1\xi_1 + \alpha_2\xi_2$ and $L_2 = \beta_1\xi_1 +
\beta_2\xi_2$ are independent, then $\mu_j\in
\Gamma(X)*I(X)$, $j=1, 2$.}

{\bf  Proof of Theorem 2.}   Necessity. Reasoning as in the proof of necessity in Theorem 1 we see that Theorem 2 fails if
  $\alpha$ does not satisfy condition  ({\ref{d1}).

Sufficiency. We follow the scheme of the proof of Theorem 1.
Note that if the torsion part of a mixed group is a bounded group,  then
  the torsion part is a direct factor of the group
    (\cite[Theorem 100.1]{Fu2}).   Taking into account this result, it follows from Lemma 6
that the group  $X$  is topologically isomorphic to a  group
of the form $\mathbb{R}^n\times D\times F$, where $n\ge 0$, and  $D$ is a torsion
free discrete Abelian group.
Denote by   $Y$  the character group of the group  $X$ and by
 $C$ the character group of the group $D$. Then $C$ is a connected compact Abelian group
(\cite[(23.17), (24.25)]{Hewitt-Ross}),  and the group  $Y$ is topologically isomorphic to the  group
$\mathbb{R}^n\times C\times H$, where $H$ is the character group of the group $F$.   In order not to complicate the
notation  we assume that
$X=\mathbb{R}^n\times D\times  F$ and  $Y=\mathbb{R}^n\times C\times H$.

 Consider the factor-group $X/(\mathbb{R}^n\times F)$, and denote by
   $[x]$ its elements.
 Since  $\mathbb{R}^n$ is a connected component of zero of the group    $X$,
 we have $\alpha(\mathbb{R}^n)=\mathbb{R}^n$.
  Since  $F$ is the subgroup   of all compact elements of   $X$, we also have
$\alpha(F)=F$. Hence,
\begin{equation}\label{Ange2}
\alpha(\mathbb{R}^n\times F)=\mathbb{R}^n\times F,
\end{equation}
and $\alpha$  induces  an  automorphism
 $\hat\alpha$  on the factor-group $X/(\mathbb{R}^n\times F)$
by the formula $\hat\alpha [x]=[\alpha x]$.
Verify that
\begin{equation}\label{1*}
{\rm Ker}(I+\hat\alpha)=\{0\}.
\end{equation}
Take $x_0\in X$ such that $[x_0]\in {\rm Ker}(I+\hat\alpha)$. We have
$(I+\hat\alpha)[x_0]=0$, i.e. $[(I+\alpha)x_0]=0$. This implies that
\begin{equation}\label{Angers1}
(I+\alpha)x_0=x'
\end{equation}
 for some   $x'\in \mathbb{R}^n\times F$. It is obvious that
$(I+\alpha)(\mathbb{R}^n\times F)\subset \mathbb{R}^n\times F$.
  It follows from $(\ref{d1})$ that the restriction of the continuous endomorphism
       $I+\alpha$ of the group $X$ to the subgroup  $\mathbb{R}^n\times F$ is a  topological
   automorphism of the group $\mathbb{R}^n\times F$. Hence,
  $x'=(I+\alpha) \tilde x$
   for some   $\tilde x\in \mathbb{R}^n\times F$. In view of  (\ref{Angers1}) this implies
     that
$(I+\alpha)x_0=(I+\alpha)\tilde x $. Taking into account
$(\ref{d1})$, it follows from this that  $x_0=\tilde x$.
 Hence,
  $[ x_0]=0$.  So, we proved
$(\ref{1*})$.

 Since $\mathbb{R}^n\times C$  is a connected component of zero of the group   $Y$,
we have  $\tilde\alpha(\mathbb{R}^n\times C)=\mathbb{R}^n\times C$.
  Since $C\times H$ is the subgroup   of all compact elements of the group    $Y$,
 we also have  $\tilde\alpha(C\times H)=C\times H$. This implies that $\tilde\alpha(C)=C$.

 Take into account that the character group of the factor-group $X/(\mathbb{R}^n\times F)$
 is topologically isomorphic in a natural way to the group $A(Y, \mathbb{R}^n\times F)=C$ (\cite[(23.25)]{Hewitt-Ross}), and note that the automorphism adjoint to  $\hat\alpha$ acts on $C$
  as
  the restriction of  $\tilde\alpha$ to $C$.
  It follows from (\ref{1*}) that the subgroup   $(I+\tilde\alpha)(C)$
 is dense in $C$ (\cite[(24.41)]{Hewitt-Ross}).  Since
   $C$ is a compact group,
$(I+\tilde\alpha)(C)=C,$ i.e. the restriction of
the topological automorphism $\tilde\alpha$ to $C$ satisfies condition (\ref{16n}).
 By Lemma 1, the symmetry of the conditional distribution of the linear form
     $L_2$ given   $L_1$ implies that the characteristic functions
$\hat\mu_j(y)$  satisfy equation $(\ref{42})$. In view of  $\tilde\alpha(C)=C$,  we can
consider the restriction of equation  $(\ref{42})$  for the characteristic functions
  $\hat \mu_j(y)$ to  $C$.
     Apply Corollary 2, considering the group $C$ instead of $Y$.
We get
\begin{equation}\label{04_04_1}
\hat\mu_j(y)=(x_j, y), \quad  y\in C,
\end{equation}
where $x_j\in D$,  $j=1, 2$.
 Substituting these expressions for $\hat\mu_j(y)$ to equation (\ref{42}) we get
 $2(x_1+\alpha x_2)\in A(X, C)=\mathbb{R}^n\times F$. This implies that
\begin{equation}\label{Z1}
x_1+\alpha x_2\in \mathbb{R}^n\times F.
\end{equation}

Consider the independent random variables $\xi'_1=\xi_1+\alpha x_2$ and
$\xi'_2=\xi_2-x_2$  with values in the
group $X$. Denote by
  $\mu'_j$   the distribution of the random variable  $\xi'_j$, $j=1, 2$.
  It is obvious that the characteristic functions    $\hat\mu'_j(y)$
  also satisfy equation $(\ref{42})$, and hence, by Lemma
   1, the conditional distribution of the linear form  $L'_2 = \xi'_1 +
\alpha\xi'_2$ given  $L'_1 = \xi'_1 +
\xi'_2$ is symmetric. It follows from (\ref{04_04_1}) and (\ref{Z1}) that  $\hat\mu'_j(y)=1$, $y\in C$, $j=1, 2$.
By Lemma 2, this implies that
$\sigma(\mu'_j)\subset A(X, C)=\mathbb{R}^n\times F$, $j=1, 2$.

  Taking into account (\ref{Ange2}), we can apply Lemma  5 to the group
  $\mathbb{R}^n\times F$ and reduce the proof of Theorem 2 to the description of distributions
     of independent random variables $\xi'_j$ with values in the group $\mathbb{R}^n\times F$,
     under assumption that the linear forms of these random variables are independent.
 Since the topological automorphism
  $\alpha$  satisfies condition ({\ref{d1}),
and $F$ is a finite Abelian group  containing no elements of order 2,
 the coefficients of the obtained linear forms, as easy to see,  are topological automorphisms of the group
       $\mathbb{R}^n\times F$. The statement of Theorem 2 follows from Lemma 7.
       $\Box$

{\bf Remark 2.}   We will verify that Theorem 2, generally speaking,
 fails if $F$ is not a finite group.
 Let  $p$ be a prime number. Denote by  $\Delta_p$ the group of $p$-adic integers,
 $\Delta_p=\{x=(x_0, x_1, \dots, x_n, \dots)\}$, where $x_j\in \{0, 1, 2, \dots, p-1\}$, $j=0, 1, 2, \dots$  (\cite[\S 10]{Hewitt-Ross}).
It is well known that  $\Delta_p$ is a torsion free totally disconnected   compact Abelian group.
 There exists a multiplication in $\Delta_p$, and       $\Delta_p$  is an integral domain with unity
 $e=(1, 0, \dots, 0,\dots)$. Note that $-e=(p-1, p-1, \dots, p-1,\dots)$. Denote by
 $\Delta_p^0={\{(x_0, x_1 \dots, x_n, \dots)\in \Delta_p: x_0\ne 0\}}$
 the set of all invertible in
  $\Delta_p$ elements.    Each automorphism       $\alpha \in {\rm
Aut}(\Delta_p)$  is of the form $\alpha x=c_\alpha x$, where
 $c_\alpha \in \Delta_p^0$, i.e.
  $\alpha$ acts as a multiplication by
$c_\alpha$  (\cite[\S 26]{Hewitt-Ross}). It follows from has been
said above that condition
({\ref{d1}) holds if and only if $\alpha\ne -I$, i.e.
$$
c_\alpha\ne(p-1, p-1, \dots, p-1, \dots).
$$

The following theorem has been proved in \cite{F}.

{\bf Theorem C}. {\it Consider the group $\Delta_p$. Let $\alpha \in {\rm
Aut}(\Delta_p)$ and $\alpha$ acts as a multiplication by $c_\alpha=(c_0, c_1, \dots, c_n, \dots)$, $c_0\ne 0$.
Then the following statements hold.

$1$. Let $\xi_1$ and  $\xi_2$ be independent random variables with values in
      the group $\Delta_p$  and distributions
  $\mu_1$ and $\mu_2$. Assume that
 the conditional distribution of the linear form  $L_2 = \xi_1 + \alpha\xi_2$
given  $L_1 = \xi_1 + \xi_2$ is symmetric. Then

$1(i)$. If $p>2$ and $c_0\ne p-1$, then
$\mu_j=m_K*E_{x_j}$, where $K$ is a compact subgroup of   $\Delta_p$, and
$x_j\in \Delta_p$, $j=1, 2$.

$1(ii)$. If $p=2$ and $c_1=0$, then $\mu_1$ and $\mu_2$ are degenerate distributions.

$2$. If one of the following conditions holds:

$2(i)$. $p>2$, $c_0=p-1$.

$2(ii)$. $p=2$, $c_1=1$,

\noindent then there exist independent random variables
$\xi_1$ and  $\xi_2$ with values in
      the group $\Delta_p$  and distributions
  $\mu_1$ and $\mu_2$ such that the conditional distribution of the
  linear form  $L_2 = \xi_1 + \alpha\xi_2$
given  $L_1 = \xi_1 + \xi_2$ is symmetric whereas
$\mu_j\notin I(\Delta_p)$, $j=1, 2$.}

Since $\Delta_p$ is a  totally disconnected   Abelian group, Gaussian distributions on
$\Delta_p$ are generated (\cite[Chapter IV, \S 6]{Pa}). Thus  theorem C  gives an
  example of a locally
compact Abelian group $X$ containing no elements of order 2 for which condition ({\ref{d1})
is necessary but not sufficient in order that the symmetry of the conditional distribution of
  the liner form $L_2=\xi_1 + \alpha \xi_2$  given  $L_1=\xi_1+\xi_2$
imply that $\mu_1$ and $\mu_2$
be convolutions
of Gaussian distributions and idempotent distributions.
Since $\Delta_p$  is
 a compact group, it coincides with the subgroup $F$   of all its compact elements.

{\bf Remark 3.}    As we noted in the proof of Theorem 2, a locally compact Abelian group $X$
such that the subgroup   of all its compact elements is finite, is topologically isomorphic to
a group of the form
\begin{equation}\label{27a}
\mathbb{R}^n\times D\times F,
\end{equation}
 where $n\ge 0$, $D$ is a torsion free discrete Abelian group, $F$ is a finite Abelian group,
 and Theorem 2 holds for groups of the form   ({\ref{27a}) under assumption that $F$ contains no elements
 of order 2, and the topological automorphism $\alpha$ satisfies condition
 ({\ref{d1}).
It follows from the proof of Theorem 2 that
if $\alpha$ satisfies condition
 ({\ref{d1}),
the proof of Theorem 2 remains  unchanged
for groups topologically isomorphic to  a group of the form $\mathbb{R}^n\times D\times F$,
where $n\ge 0$, $D$ is a torsion free discrete Abelian group,   $F$ is a torsion discrete Abelian group such that $F$ contains no elements
 of order 2  and has the following property: each endomorphism of $F$
 with zero kernel is an automorphism.
The only difference is at the end of the proof, where instead of Lemma 7
 we should use the same result for the groups of the form
 $X=\mathbb{R}^n\times G$, where $G$ is an arbitrary discrete Abelian
 group (\cite[Remark 13.19]{Fe5}).

 As has been observed in  \cite[Remark 17.27]{Fe5}, if we assume that a topological automorphism
 $\alpha$ of a group $X$  satisfies  stronger  conditions than
  ({\ref{d1}), namely $(\ref{E})$, then
 the statement of Theorem  2 holds true for a much wider class   of group    $X$, than
 groups of the form ({\ref{27a}).  They are groups topologically isomorphic to a group
  of the form
 \begin{equation}\label{XXX}
\mathbb{R}^n\times G,
\end{equation}
 where $n\ge 0$,  and $G$ is an arbitrary discrete Abelian group containing no elements of order 2.

 It would be interesting to find out  if the statement of Theorem 2 is valid for the groups of the form
 ({\ref{XXX}),  if we change of condition ({\ref{E}) for
   ({\ref{d1}).

\bigskip

\centerline{\textbf{4.  One more characterization of the Gaussian }}
\centerline{\textbf{distribution on locally compact Abelian groups}}

\bigskip

The following characterization theorem was proved in
\cite{Fe4}, see also \cite[Theorem 16.2]{Fe5}.

{\bf Theorem D}. {\it Let  $X$ be a locally compact Abelian group
containing
no elements of order $2$.
Let
$\alpha_j$, $\beta_j$ be topological automorphisms of
the group   $X$
such that
\begin{equation}\label{1}
\beta_i\alpha_i^{-1} \pm \beta_j\alpha_j^{-1} \in {\rm
Aut}(X)
\end{equation}
  for all $i \ne j$.
Let
$\xi_j$, $j = 1, 2,\dots, n$, $n \ge 2,$ be independent random
variables with values in $X$ and distributions  $\mu_j$
with non-vanishing characteristic functions.  If the conditional distribution of
  the linear form   $L_2 = \beta_1\xi_1 + \dots + \beta_n\xi_n$
given $L_1 = \alpha_1\xi_1 + \dots + \alpha_n\xi_n$
is symmetric, then all $\mu_j$ are Gaussian distributions.}

It is easy to see that if $n=2$, it suffices to prove
 Theorem D in the case
 when $L_1 = \xi_1 +\xi_2$,
$L_2 = \xi_1 +\alpha\xi_2$, where $\alpha\in {\rm Aut}(X)$.
Condition   (\ref{1}) in this case is transformed into condition ({\ref{E}).
 We prove in this section that if $n=2$, then Theorem D
  holds under a much weaker condition than
   (\ref{E}), namely (\ref{d1}).

 {\bf Theorem 3.}
 {\it  Let  $X$ be a locally compact Abelian group containing
no elements of order $2$.  Let $\alpha$ be a topological automorphism of the group
$X$.
Let  $\xi_1$ and  $\xi_2$ be independent random variables with values in
$X$   and distributions   $\mu_1$ and $\mu_2$ with non-vanishing characteristic functions.
The symmetry of the conditional distribution of the linear form
$L_2 = \xi_1 + \alpha\xi_2$ given  $L_1 = \xi_1 +
\xi_2$ implies that   $\mu_1$ and $\mu_2$ are Gaussian distributions if and only if
  $\alpha$ satisfies condition  $(\ref{d1})$.}

To prove Theorem 3 we need two lemmas.

The following lemma is standard and it was proved
in (\cite[Lemma 10.1]{Fe5}) under assumption that
$\alpha_j$, $\beta_j$ are topological automorphisms. Its proof remains valid
if $\alpha_j$, $\beta_j$ are continuous endmorphisms. We formulate the lemma for two independent random variables.

{\bf Lemma 8}. {\it Let $X$ be a locally compact Abelian group, $Y$ be
its character group. Let
$\alpha_j$, $\beta_j$ be continuous endomorphisms of the group $X$. Let $\xi_1$ and $\xi_2$ be independent random variables with values in the group
       $X$  and distributions
 $\mu_1$ and $\mu_2$.
 The linear forms $L_1 = \alpha_1\xi_1 +
\alpha_2\xi_2$ and $L_2 = \beta_1\xi_1  + \beta_2\xi_2$ are independent
if and only if the characteristic functions $\hat\mu_j(y)$ satisfy the
equation}
\begin{equation}\label{4}
\hat\mu_1(\tilde\alpha_1 u+\tilde\beta_1 v)\hat\mu_2(\tilde\alpha_2 u+\tilde\beta_2 v)=\hat\mu_1(\tilde\alpha_1
u)\hat\mu_2(\tilde\alpha_2
u)\hat\mu_1(\tilde\beta_1 v)\hat\mu_2(\tilde\beta_2 v),
\quad u, v \in
Y.
\end{equation}

The following statement is a group analog of the well-known Cramer
theorem about decomposition of the Gaussian distribution.

 {\bf Lemma 9} (\cite{Fe1}, see also \cite[Theorem 4.6]{Fe5}). {\it Let $X$ be a locally compact  Abelian containing no subgroup topologically isomorphic
to the circle group $\mathbb{T}.$  Let
 $\mu \in \Gamma(X)$. If
  $\mu=\lambda_1*\lambda_2$,
 where $\lambda_j\in {\rm M}^1(X)$,
 then $\lambda_j \in \Gamma(X)$, $j=1, 2$.}

  {\bf Proof of Theorem 3.}  The proof of necessity is the same as in Theorem 1.

 Sufficiency. Denote by   $Y$  the character group of the group  $X$.
By Lemma 5,  the symmetry of the conditional distribution
of the linear form $L_2$
given $L_1$ implies that the linear forms
$M_1=(I+\alpha)\xi_1+2\alpha\xi_2$ and
$M_2=2\xi_1+(I+\alpha)\xi_2$ are independent.
By Lemma 8, this implies that the characteristic functions
$\hat\mu_j(y)$ satisfy equation $(\ref{4})$,
which takes the form
\begin{equation}\label{5}
\hat\mu_1((I+\tilde\alpha) u+2 v)\hat\mu_2(2\tilde\alpha  u+(I+\tilde\alpha) v)=\hat\mu_1((I+\tilde\alpha)
u)\hat\mu_2(2\tilde\alpha
u)\hat\mu_1(2 v)\hat\mu_2((I+\tilde\alpha) v), \quad u, v \in
Y.
\end{equation}
Set $\nu_j = \mu_j* \bar \mu_j$. Then  $\hat \nu_j(y) = |\hat \mu_j(y)|^2 > 0,$  $y \in Y$, $j=1, 2$.
 It is obvious that the characteristic functions
 $\hat \nu_j(y)$
also satisfy equation $(\ref{5})$.

Put $\psi_j(y)=-\log\hat\nu_j(y)$, $j=1, 2$. It follows from
$(\ref{5})$ that the functions $\psi_j(y)$ satisfy the equation
\begin{equation}\label{6}
\psi_1((I+\tilde\alpha) u+2 v)+\psi_2(2\tilde\alpha  u+(I+\tilde\alpha) v)=P(u)+Q(v), \quad u, v \in
Y,
\end{equation}
where  \begin{equation}\label{16.09.15.1}
P(y)=\psi_1((I+\tilde\alpha)
y)+\psi_2(2\tilde\alpha
y), \quad Q(y)=\psi_1(2 y)+\psi_2((I+\tilde\alpha) y).
\end{equation}

We use the finite difference method to solve equation (\ref{6}) (see e.g. \cite[\S 10]{Fe5}).
Let $h_1$ be an arbitrary element of the group $Y$.
Substitute in
 (\ref{6}) $u+(I+\tilde\alpha) h_1$ for $u$ and $v-2\tilde\alpha  h_1$ for $v$. Subtracting equation
  (\ref{6}) from the resulting equation, we get
  \begin{equation}\label{7}
    \Delta_{(I-\tilde\alpha)^2 h_1}{\psi_1((I+\tilde\alpha) u+2 v)}
    =\Delta_{(I+\tilde\alpha) h_1} P(u)+\Delta_{-2\tilde\alpha  h_1} Q(v),
\quad u, v\in Y.
\end{equation}
Let $h_2$ be an arbitrary element of the group $Y$.
Substitute in
(\ref{7})   $u+2h_{2}$ for $u$ and
$v-(I+\tilde\alpha)h_{2}$ for $v$. Subtracting equation
  (\ref{7}) from the resulting equation, we obtain
 \begin{equation}\label{8}
     \Delta_{2 h_2}\Delta_{(I+\tilde\alpha) h_1} P(u)+\Delta_{-(I+\tilde\alpha) h_2}\Delta_{-2\tilde\alpha  h_1} Q(v)=0,
\quad u, v\in Y.
\end{equation}
Let $h$ be an arbitrary element of the group $Y$.
Substitute in
(\ref{8})   $u+h$ for $u$. Subtracting equation
  (\ref{8}) from the resulting equation, we get
 \begin{equation}\label{9}
   \Delta_{h}\Delta_{2 h_2}\Delta_{(I+\tilde\alpha) h_1} P(u)=0,
\quad u\in Y.
\end{equation}

It is worth noting that we got (\ref{9})  using neither $X$ contains no elements of
order 2, nor $\alpha$ satisfies condition (\ref{d1}).

Since the group $X$ contains no elements of order 2,
the subgroup $Y^{(2)}$ is dense in $Y$ (\cite[$($24.22)]{Hewitt-Ross}). Since
$\alpha$ satisfies condition
$(\ref{d1})$, the subgroup
$(I+\tilde\alpha)(Y)$ is also dense in
$Y$ (\cite[(24.41)]{Hewitt-Ross}). Taking into account that
$h, h_1, h_2$ are arbitrary elements of $Y$, it follows from
(\ref{9}) that the function $P(y)$ satisfies the equation
\begin{equation}\label{10}
\Delta_h^{3} P(y) = 0, \quad y, h  \in Y.
\end{equation}
Taking into account that $P(-y)=P(y)$, $y\in Y$ and $P(0)=0$, it follows from
(\ref{10}) that
\begin{equation}\label{18.12.1}
P(u+v)+P(u-v)=2[P(u)+P(v)], \quad u, v\in Y.
\end{equation}

Let $\mu$ be the distribution of the random
variable $M_1=(I+\alpha)\xi_1+2\alpha\xi_2$. Put
$\nu=\mu*\bar\mu$. It is obvious that  $\mu=(I+\alpha)(\mu_1)*(2\alpha)(\mu_2)$. It follows from
(\ref{16.09.15.1}) that
 the characteristic function $\hat\nu(y)$ is of the form $$\hat\nu(y)=\exp\{-P(y)\}, \quad y\in Y.$$
In view of (\ref{18.12.1}), we have $\nu\in \Gamma(X)$.
Since the group $X$ contains no elements of
order 2,  $X$ contains no subgroup topologically isomorphic to
the circle group ${\mathbb T}$.
By Lemma 9,
we get
  that
 $(I+\alpha)(\mu_1)\in \Gamma(X)$. Since, by condition
$(\ref{d1})$,
 $I+\alpha$ is a continuous monomorphism itself,
  $\mu_1\in \Gamma(X)$.

 Let $k$ be an arbitrary element of the group $Y$.
Substitute in
(\ref{8})   $v+k$ for $u$. Subtracting equation
  (\ref{8}) from the resulting equation, we get
  \begin{equation}\label{9a}
   \Delta_{k}\Delta_{-(I+\tilde\alpha) h_2}\Delta_{-2\tilde\alpha  h_1} Q(v)=0,
\quad v\in Y.
\end{equation}
Taking into account (\ref{16.09.15.1}) and $(\ref{9a})$, considering the distribution
of the random variable $M_2=2\xi_1+(I+\alpha)\xi_2$ and
reasoning similarly as above, we obtain that
$\mu_2\in \Gamma(X)$.  $\Box$

\vskip 1 cm

\noindent G. M. Feldman

\medskip

\noindent B. Verkin Institute for Low Temperature \\
Physics and Engineering of the\\
National Academy of Sciences of Ukraine\\
47, Lenin Ave  \\
Kharkiv, 61103  \\
Ukraine

\end{document}